\documentclass[12pt]{article}
\usepackage{mathrsfs}

\usepackage{t1enc}
\usepackage[latin1]{inputenc}
\usepackage[english]{babel}
\usepackage{amsmath,amsthm,amssymb}
\usepackage{amsfonts}
\usepackage{latexsym}
\usepackage{enumerate}

\usepackage{graphicx}
\usepackage[natural]{xcolor}
\usepackage{tikz}
\usepackage{tikz-3dplot}
\usetikzlibrary{shapes}
\usetikzlibrary{arrows,decorations.pathmorphing,backgrounds,positioning,fit,petri,automata}
\usetikzlibrary {positioning}

\theoremstyle{plain}
\newtheorem{thm}{Theorem}[section]

\newtheorem{lem}[thm]{Lemma}

\newtheorem{claim}[thm]{Claim}

\theoremstyle{plain}

\theoremstyle{plain}

\theoremstyle{plain}

\title{Erd\H{o}s-Gy\'{a}rf\'{a}s Conjecture for $P_8$-free graphs}
\author{Yuping Gao$^{a}$, Songling Shan$^{b}$\\
{\small a. School of Mathematics and Statistics, Lanzhou University, Lanzhou 730000, China}\\
{\small b. Department of Mathematics, Illinois State University, Normal, IL 61790, USA}}
\date{}

\begin{document}
\baselineskip 0.65cm

\maketitle
\begin{abstract} A graph  is $P_8$-free if it contains no induced subgraph  isomorphic to the path $P_8$ on eight vertices. In 1995,  Erd\H{o}s and  Gy\'{a}rf\'{a}s conjectured that every graph of minimum degree at least three contains a cycle whose length is a power of two. In this paper, we confirm the conjecture for $P_8$-free graphs by showing that there exists a cycle of length four or eight in every $P_8$-free graph with minimum degree at least three.

\medskip

\end{abstract}

\section{Introduction}
All graphs considered in this paper are undirected and simple. Let $G$ be a graph. The vertex set, the edge set, the maximum degree and the minimum degree of $G$ are denoted by $V(G)$, $E(G)$, $\Delta(G)$ and $\delta(G)$, respectively. For a vertex $v\in V(G)$, the set of neighbors of $v$ in $G$ is denoted by $N_{G}(v)$ or $N(v)$ if $G$ is understood.
Let $S\subseteq V(G)$, we use $G[S]$ to denote the subgraph of $G$ induced by $S$ and
$G-S$ to denote the subgraph $G[V(G)\setminus S]$.  We write $u\thicksim v$  if $uv\in E(G)$ and $u\nsim v$ otherwise. The connectivity of $G$ is denoted by $\kappa(G)$. A $uv$-\emph{path} is a path having ends as $u$ and $v$. Let $P$ be a path and $x,y\in V(P)$, we use $xPy$ to denote the subpath of $P$ with ends $x$ and $y$.

 A path on  $k$ vertices is denoted by $P_k$. A cycle on $k$ vertices
 is denoted by $C_k$ and is called a $k$-\emph{cycle}.
 The length of a path or cycle is the number of edges it contains.
 The well-known Erd\H{o}s-Gy\'{a}rf\'{a}s Conjecture~\cite{E1995} states that every graph of minimum degree at least three contains a $2^{m}$-cycle for some integer $m\geq 2$. The conjecture is confirmed for some graph classes including $K_{1,m}$-free graphs of  minimum degree at least $m+1$ or maximum degree at least $2m-1$~\cite{S1998}, 3-connected cubic planar graphs~\cite{HK2013},  planar claw-free graphs~\cite{DS2001} and some Cayley graphs~\cite{GM2018,GV2021}. In~\cite{NESB2014}, it is proved that every cubic claw-free graph  contains a cycle whose length is $2^{k}$, or $3\cdot 2^{k}$, for some positive integer $k$.

Given a graph $H$, a graph $G$ is  $H$-\emph{free} if $G$ does not contain any induced subgraph  isomorphic to $H$. In this paper, we confirm Erd\H{o}s-Gy\'{a}rf\'{a}s Conjecture for $P_8$-free graphs by showing the following two theorems.

\begin{thm}\label{thm} Every  $P_5$-free graph  with minimum degree at least three contains a $4$-cycle.
\end{thm}

\begin{thm}\label{thm2} Every  $P_8$-free graph with minimum degree at least three contains a $4$-cycle or $8$-cycle.
\end{thm}

In confirming the Erd\H{o}s-Gy\'{a}rf\'{a}s Conjecture for $P_8$-free graphs, Theorem~\ref{thm2} alone suffices. But
 we include Theorem~\ref{thm} as it is stronger than the restriction of Theorem~\ref{thm2} on $P_5$-free graphs and also its
proof technique may be of independent interests.

The remainder of the paper is organized as follows. In Section~2, we prove Theorem~\ref{thm}. In Section~3, we prove Theorem~\ref{thm2}.

\section{Proof of Theorem~\ref{thm}}

\begin{proof}[Proof of Theorem~\rm{\ref{thm}}]  Let $G$ be a $P_5$-free graph with $\delta(G) \ge 3$.
	We may assume that $G$ is connected. Otherwise, we  consider a component  of $G$ instead. Furthermore, assume that $G$ is not complete and  $G$ contains no $C_4$ since otherwise we are done. Let $S$ be a minimum cut-set of $G$.
For $x\in S$, 	
a component $D$ of $G-S$ is a \emph{complete neighborhood component} (\emph{CNC}) of $x$ if $x$ is adjacent in $G$ to every vertex of $D$, and $D$
is a  \emph{non-CNC} of $x$ otherwise.
We need the following claim.

\begin{claim}\label{claim1} {\rm(i)} For any vertex $x\in S$ and any component $D$ of $G-S$, $N_G(x)\cap V(D) \ne \emptyset$.

{\rm(ii)} For any vertex $x\in S$, $x$ has at least $c(G-S)-1$ CNCs. Equivalently, $x$ has at most one non-CNC.

{\rm(iii)} For any CNC $D$ of a vertex $x\in S$, $|V(D)|\leq 2$.

{\rm(iv)} $|S|\geq 2$, i.e., $\kappa(G)\geq 2$.

Moreover, let $x,y\in S$ be distinct. Then the following statements hold.

{\rm(v)} If $x$ and $y$ have a common CNC $D$, then $|V(D)|=1$.

{\rm(vi)} $x$ and $y$ have at most one common CNC.

{\rm(vii)} $c(G-S)\leq 3$. Furthermore, if $c(G-S)=3$, then $x$ and $y$ have exactly one common CNC.

{\rm(viii)} There exist two  vertices from $S$ that are nonadjacent in $G$.

\end{claim}

\begin{proof}
For Statement (i), as $S$ is a minimum cut-set of $G$, for any $x\in S$ and any component $D$ of $G-S$,  it follows that $N_G(x)\cap V(D) \ne \emptyset$.
	
For Statement (ii),  suppose  instead that $x$ has two non-CNCs $D_1$ and $D_2$. Let $u_i\in V(D_i)$ such that $x\nsim u_i$, and $P_i$ be a shortest
	path of $D_i$ from $u_i$ to a neighbor, say $x_i$ of $x$ in $G$
	from $V(D_i)$, $i=1,2$. By the choice,  $P_i$  is an induced path of $D_i$
	such that the only vertex of $P_i$ that is adjacent in $G$ to $x$
	is $x_i$, $i=1,2$. 	Then $u_1P_1x_1xx_2P_2u_2$ contains an induced $P_5$, contradicting $G$
	being $P_5$-free.

For Statement	(iii),  suppose  instead that $|V(D)|\geq 3$. If $D$ is a complete graph, then $G[V(D)\cup \{x\}]$ contains a $C_4$ and we are done. So assume that $D$ is not a complete graph. Then $D$ contains an induced $P_3$. Since $x$ is adjacent to every vertex in $D$, especially adjacent to every vertex in the $P_3$. It follows that $G$ contains a $C_4$.

For Statement (iv),  by (ii), each vertex $x\in S$ has a CNC $D$. By (iii), $|V(D)|\leq 2$. Let $u\in V(D)$. Then $u$ has a neighbor $y\in S\setminus \{x\}$ by $\delta(G)\geq 3$. So $|S|\geq 2$.

Statements (v) and (vi) follow by the assumption that $G$ contains no $C_4$.

For Statement (vii),  by (ii), each of $x$ and $y$ has at most one non-CNC. Then $x$ and $y$ have at least $c(G-S)-2$ common CNCs. This is a contradiction to statement (vi) if $c(G-S)\geq 4$. Furthermore, if $c(G-S)=3$, then $x$ and $y$ have exactly one common CNC.

For Statement (viii),  by  (ii) and (iii), $x$ has a CNC $D$ with $|V(D)|\leq 2$. Let $u\in V(D)$.  If $u$ has at least three neighbors from $S$, then there exist two nonadjacent vertices in $N(u)\cap S$ since $G$ contains no $C_4$. So $u$ has at most two neighbors in $S$.
This, together with $\delta(G) \ge 3$, implies  $|V(D)|=2$. Let $v$ be the neighbor of $u$
in $D$.
Since $|V(D)|\leq 2$
and $\delta(G) \ge 3$, $u$ has a neighbor $ y\in S\setminus \{x\}$.
Furthermore, $y\nsim x$, for otherwise $yuvxy$ is a $C_4$ of $G$.
%
\end{proof}

By Claim~\ref{claim1}(viii), we let $x,y\in S$ such that $x\nsim y$ in $G$.

\begin{claim}\label{claim2} $c(G-S)=2$.
\end{claim}

\begin{proof} By Claim~\ref{claim1}(vii), $c(G-S)\leq 3$. Suppose that $c(G-S)=3$. By Claim~\ref{claim1}(vii), $x$ and $y$ have exactly one common CNC $D_1$. By Claim~\ref{claim1}(ii), let $D_2$ be a CNC of $x$, $D_3$ be a CNC of $y$. Note that $D_2\neq D_3$. Let $u_i\in V(D_i),i=1,2,3$. Then $u_3\nsim x, u_2\nsim y$ since $G$ contains no $C_4$. It follows that $u_3yu_1xu_2$ is an induced $P_5$ in $G$, a contradiction.
 \end{proof}

By Claim~\ref{claim2}, $c(G-S)=2$. Let $D_1,D_2$ be the two components of $G-S$.
 Assume first that $x$ and $y$ have no common CNC. By Claim~\ref{claim1}(ii), assume by symmetry that $D_1$ is a CNC of $x$ and $D_2$ is a CNC of $y$.  By Claim~\ref{claim1}(i), there exist $y_1\in N(y)\cap V(D_1)$ and $x_1\in N(x)\cap V(D_2)$.
 Then $xy_1yx_1x$ is a $C_4$, giving a contradiction.

Now assume that $x$ and $y$ have a common CNC $D_1$ and let $u\in V(D_1)$. Then $|V(D_1)|=1$ by Claim~\ref{claim1}(v). By Claim~\ref{claim1}(i), there exist $x_1\in N(x)\cap V(D_2)$ and $y_1\in N(y)\cap V(D_2)$. Note that $x_1\ne y_1$. If $x_1\nsim y_1$, then $x_1xuyy_1$ is an induced $P_5$ in $G$. So $x_1\thicksim y_1$. Since $\delta(G)\geq 3$, $u$ has a neighbor $z\in S\setminus \{x,y\}$. Because $G$ has no $C_4$, we have $z\nsim x_1,y_1$.
Now it must be the case that $z\thicksim x$, as otherwise $zuxx_1y_1$ is an induced $P_5$ in $G$. Similarly,  $z\sim y$. However, it follows that $z$ is a common neighbor of $x$
and $y$ other than $u$, showing a contradiction. We complete the proof of Theorem~\ref{thm}.
\end{proof}

\section{Proof of Theorem~\ref{thm2}}

The Lemma below was shown in \cite{NESB2014}.

\begin{lem}\label{lem2.1} Let $G$ be a graph with $\delta(G)\geq 3$. If $G$ does not contain $C_4$, then  $G$ has an induced cycle $C_k$ for some $k\geq 5$.
\end{lem}

\begin{proof}[Proof of Theorem~\rm{\ref{thm2}}] Let $G$ be a $P_8$-free graph with $\delta(G) \ge 3$.
	We may assume that $G$ is connected. Otherwise, we just consider a component  of $G$. Furthermore, assume that $G$ contains neither $C_4$ nor $C_8$ since otherwise we are done. By Lemma~\ref{lem2.1}, $G$ contains an induced $C_k$ for some $k\geq 5$. Let $C=v_1v_2\ldots v_k v_1$ be a shortest induced cycle in $G$ of length at least 5. Then $5\leq k\leq 7$ since $G$ is $P_8$-free and $G$ contains neither $C_4$ nor $C_8$.

\begin{claim}\label{claim2.1} If $k=5$, then no two consecutive vertices on $C$ share a common neighbor in $G$.
\end{claim}

\begin{proof} Suppose the claim does not hold. We assume, without loss of generality, that $v_1$ and $v_2$ have a common neighbor $v_6\not\in\{v_1,v_2,\ldots,v_5\}$. Then $v_6\nsim v_i$ for $i\in\{3,4,5\}$ as $G$ contains no $C_4$. We conclude that $v_6$ has a neighbor $v_7\not\in\{v_1,v_2,\ldots,v_6\}$ since $\delta(G)\geq 3$. The minimum degree condition is repeatedly used in the following proof and we omit the reason in the following when we say that $v_i$ has a neighbor $v_j$ for $i\neq j$. It can be seen that $v_7\nsim v_i$ for $i\in \{1,2,3,5\}$ as $G$ contains no $C_4$.

{\bf Case 1} $v_7\thicksim v_4$.

In this case, $v_7$ has a neighbor $v_8\not\in\{v_1,v_2,\ldots,v_7\}$. And $v_8\nsim v_i$ for $i\in\{1,2,3,5\}$ since $G$ contains no $C_4$ or $C_8$.

{\bf Subcase 1.1} $v_8\thicksim v_6$.

In this case, $v_8\nsim v_4$ otherwise $v_8v_4v_7v_6v_8$ is a $C_4$. It follows that $v_8$ has a neighbor $v_9\not\in \{v_1,v_2,\ldots,v_8\}$. And $v_9\nsim v_1$ otherwise $v_9v_1v_6v_8v_9$ is a $C_4$, $v_9\nsim v_2$ otherwise $v_9v_2v_6v_8v_9$ is a $C_4$, $v_9\nsim v_3$ otherwise $v_9v_8v_7v_4v_5v_1v_2v_3v_9$ is a $C_8$, $v_9\nsim v_4$ otherwise $v_9v_8v_7v_4v_9$ is a $C_4$, $v_9\nsim v_5$ otherwise $v_9v_8v_7v_6v_2v_3v_4v_5v_9$ is a $C_8$,  $v_9\nsim v_6$ otherwise $v_9v_8v_7v_6v_9$ is a $C_4$,  $v_9\nsim v_7$ otherwise $v_9v_8v_6v_7v_9$ is a $C_4$. So $v_9$ has two neighbors $v_{10},v_{11}\not\in\{v_1,v_2,\cdots,v_{9}\}$. At least one of $v_{10}$ and $v_{11}$, say $v_{10}$, is not adjacent to $v_8$. Note that $v_{10}\nsim v_7$ otherwise $v_{10}v_9v_8v_7v_{10}$ is a $C_4$, $v_{10}\nsim v_4$ otherwise $v_{10}v_9v_8v_7v_6v_2v_3v_4v_{10}$ is a $C_8$,  $v_{10}\nsim v_5$ otherwise $v_{10}v_9v_8v_7v_6v_2v_1v_5v_{10}$ is a $C_8$,  $v_{10}\nsim v_1$ otherwise $v_{10}v_9v_8v_7v_4v_3v_2v_1v_{10}$ is a $C_8$.  Similarly, $v_{10} \not\sim v_2$. So $v_{10}v_9v_8v_7v_4v_5v_1v_2$ is an induced $P_8$ in $G$, a contradiction. (See Figure~1(a) for an illustration.)

\begin{center}
\begin{tikzpicture}

{\tikzstyle{every node}=[draw,circle,fill=white,minimum size=4pt,
                            inner sep=0pt]
    \draw (0,6) node (v1) {}
        -- ++(0:1cm) node (v2){}
        -- ++(0:1cm) node (v3){}
       -- ++(0:1cm) node (v4){}
        -- ++(0:1cm) node (v5){};

 \draw (v1)  {}
 -- ++(300:1cm) node (v6)   {}
  -- ++(270:1cm) node (v7)   {}
 --++(270:1cm) node (v8)   {}
  -- ++(270:1cm) node (v9)   {}
  -- ++(240:1cm) node (v10)   {} ;
 \draw (v9)  {}
 -- ++(300:1cm) node (v11)   {} ;
 \draw (v2)  -- (v6) ;
 \draw[red] (v10)-- (v9) -- (v8)-- (v7)-- (v4)-- (v5);
 \draw [color=red](v2)--(v1) .. controls (1,7) and (4,6.5)..(v5) ;
 \draw (v6) .. controls (0,4) ..(v8) ;
        }

 \draw (v1) node [left] {$v_1$};
 \draw (v2) node [below right] {$v_2$};
 \draw (v3) node [below] {$v_3$};
 \draw (v4) node [below] {$v_4$};
 \draw (v5) node [below] {$v_5$};
 \draw (v6) node [right] {$v_6$};
 \draw (v7) node [below right] {$v_7$};
 \draw (v8) node [right] {$v_8$};
 \draw (v9) node [right] {$v_9$};
 \draw (v10) node [left] {$v_{10}$};
 \draw (v11) node [right] {$v_{11}$};

\node at (1,0.5) {(a) Subcase 1.1};

{\tikzstyle{every node}=[draw,circle,fill=white,minimum size=4pt,
                            inner sep=0pt]
    \draw (6,6) node (v1) {}
        -- ++(0:1cm) node (v2){}
        -- ++(0:1cm) node (v3){}
       -- ++(0:1cm) node (v4){}
        -- ++(0:1cm) node (v5){};

 \draw (v1)  {}
 -- ++(300:1cm) node (v6)   {}
  -- ++(270:1cm) node (v7)   {}
 --++(270:1cm) node (v8)   {}
  -- ++(240:1cm) node (v9)   {}
  -- ++(240:1cm) node (v11)   {} ;
 \draw (v8)  {}
 -- ++(300:1cm) node (v10)   {} ;
 \draw (v9)  {}
 -- ++(300:1cm) node (v12)   {} ;
 \draw (v2)  -- (v6) ;
 \draw (v7)  -- (v4);
 \draw[red] (v11) --(v9) -- (v8)-- (v7)-- (v4)-- (v5);
  \draw [color=red](v2)--(v1) .. controls (7,7) and (10,6.5)..(v5) ;
        }

 \draw (v1) node [left] {$v_1$};
 \draw (v2) node [below right] {$v_2$};
 \draw (v3) node [below] {$v_3$};
 \draw (v4) node [below] {$v_4$};
 \draw (v5) node [below] {$v_5$};
 \draw (v6) node [left] {$v_6$};
 \draw (v7) node [left] {$v_7$};
 \draw (v8) node [left] {$v_8$};
 \draw (v9) node [left] {$v_9$};
 \draw (v10) node [right] {$v_{10}$};
\draw (v11) node [left] {$v_{11}$};
\draw (v12) node [right] {$v_{12}$};

\node at (7,0.5) {(b) Subcase 1.2};
\end{tikzpicture}
\end{center}

\begin{center}
\begin{tikzpicture}

{\tikzstyle{every node}=[draw,circle,fill=white,minimum size=4pt,
                            inner sep=0pt]
    \draw (0,6) node (v1) {}
        -- ++(0:1cm) node (v2){}
        -- ++(0:1cm) node (v3){}
       -- ++(0:1cm) node (v4){}
        -- ++(0:1cm) node (v5){};

 \draw (v1)  {}
 -- ++(300:1cm) node (v6)   {}
  -- ++(270:1cm) node (v7)   {}
 --++(270:1cm) node (v8)   {}
  -- ++(270:1cm) node (v9)   {}
  -- ++(240:1cm) node (v10)   {}
  -- ++(240:1cm) node (v12)   {};
 \draw (v9)  {}
 -- ++(300:1cm) node (v11)   {} ;
 \draw (v10)  {}
 -- ++(300:1cm) node (v13)   {} ;
 \draw (v2)  -- (v6) ;
 \draw (v7)  -- (v4) -- (v8) ;
 \draw (v1).. controls (1,7) and (4,6.5)..(v5);
        }
\node at (1,-0.3) {(c) Subcase 1.2};

\begin{scope}[shift={(0,-1.5)}]

 \draw (v1) node [left] {$v_1$};
 \draw (v2) node [below right] {$v_2$};
 \draw (v3) node [below] {$v_3$};
 \draw (v4) node [below] {$v_4$};
 \draw (v5) node [below] {$v_5$};
 \draw (v6) node [left] {$v_6$};
 \draw (v7) node [left] {$v_7$};
 \draw (v8) node [left] {$v_8$};
 \draw (v9) node [left] {$v_9$};
 \draw (v10) node [left] {$v_{10}$};
 \draw (v11) node [right] {$v_{11}$};
\draw (v12) node [left] {$v_{12}$};
\draw (v13) node [right] {$v_{13}$};

{\tikzstyle{every node}=[draw,circle,fill=white,minimum size=4pt,
                            inner sep=0pt]
    \draw (6,6) node (v1) {}
        -- ++(0:1cm) node (v2){}
        -- ++(0:1cm) node (v3){}
       -- ++(0:1cm) node (v4){}
        -- ++(0:1cm) node (v5){};

 \draw (v1)  {}
 -- ++(300:1cm) node (v6)   {}
  -- ++(270:1cm) node (v7)   {}
 --++(240:1cm) node (v8)   {}
  -- ++(240:1cm) node (v10)   {} ;
 \draw (v7)  {}
 -- ++(300:1cm) node (v9)   {} ;
  \draw (v8)  {}
 -- ++(300:1cm) node (v11)   {} ;
 \draw (v2)  -- (v6) ;
 \draw (v1) .. controls (7,7) and (10,6.5)..(v5) ;

        }

 \draw (v1) node [left] {$v_1$};
 \draw (v2) node [below right] {$v_2$};
 \draw (v3) node [below] {$v_3$};
 \draw (v4) node [below] {$v_4$};
 \draw (v5) node [below] {$v_5$};
 \draw (v6) node [left] {$v_6$};
 \draw (v7) node [left] {$v_7$};
 \draw (v8) node [left] {$v_8$};
 \draw (v9) node [right] {$v_9$};
 \draw (v10) node [left] {$v_{10}$};
 \draw (v11) node [right] {$v_{11}$};
\end{scope}
\node at (7,-0.3) {(d) Case 2};
\node at (5,-1.5) {{\bf Figure 1:} Illustration for Claim~\ref{claim2.1}};
\end{tikzpicture}
\end{center}

{\bf Subcase 1.2} $v_8\nsim v_6$.

If $v_8\nsim v_4$, then $v_8$ has two neighbors $v_{9},v_{10}\not\in \{v_1,v_2,\ldots,v_8\}$. At least one of $v_9$ and $v_{10}$, say $v_9$, is not adjacent to $v_7$. Moreover, $v_9\nsim v_1$ otherwise $v_9v_8v_7v_4v_3v_2v_6v_1v_9$ is a $C_8$, $v_9\nsim v_2$ otherwise $v_9v_8v_7v_4v_5v_1v_6v_2v_9$ is a $C_8$, $v_9\nsim v_3$ otherwise $v_9v_8v_7v_4v_5v_1v_2v_3v_9$ is a $C_8$, $v_9\nsim v_4$ otherwise $v_9v_8v_7v_4v_9$ is a $C_4$, $v_9\nsim v_5$ otherwise $v_9v_8v_7v_6v_2v_3v_4v_5v_9$ is a $C_8$, $v_9\nsim v_6$ otherwise $v_9v_8v_7v_6v_9$ is a $C_4$. So $v_9$ has two neighbors $v_{11},v_{12}\not\in \{v_1,v_2,\ldots,v_9\}$. At least one of $v_{11}$ and $v_{12}$, say $v_{11}$, is not adjacent to $v_8$. Moreover, $v_{11}\nsim v_7$ otherwise $v_{11}v_9v_8v_7v_{11}$ is a $C_4$, $v_{11}\nsim v_4$ otherwise $v_{11}v_9v_8v_7v_6v_2v_3v_4v_{11}$ is a $C_8$, $v_{11}\nsim v_5$ otherwise $v_{11}v_9v_8v_7v_6v_2v_1v_5v_{11}$ is a $C_8$, $v_{11}\nsim v_1$ otherwise $v_{11}v_9v_8v_7v_4v_3v_2v_1v_{11}$ is a $C_8$, $v_{11}\nsim v_2$ otherwise $v_{11}v_9v_8v_7v_4v_5v_1v_2v_{11}$ is a $C_8$. It follows that $v_{11}v_9v_8v_7v_4v_5v_1v_2$ is an induced $P_8$ in $G$, a contradiction. (See Figure~1(b) for an illustration.)

Now assume that $v_8\thicksim v_4$. Then $v_8$ has a neighbor $v_9\not\in \{v_1,v_2,\ldots,v_8\}$. Furthermore, $v_9\nsim v_i$ for $i\in\{1,2,\ldots,6\}$ same as the case when $v_8\nsim v_4$. And $v_9\nsim v_7$ otherwise $v_9v_8v_4v_7v_9$ is a $C_4$. It follows that $v_9$ has two neighbors $v_{10},v_{11}\not\in\{v_1,v_2,\ldots,v_9\}$. At least one of $v_{10}$ and $v_{11}$, say $v_{10}$, is not adjacent to $v_8$. Moreover, $v_{10}\nsim v_1$ otherwise $v_{10}v_9v_8v_7v_4v_3v_2v_1v_{10}$ is a $C_8$, $v_{10}\nsim v_2$ otherwise $v_{10}v_9v_8v_4v_7v_6v_1v_2v_{10}$ is a $C_8$, $v_{10}\nsim v_3$ otherwise $v_{10}v_9v_8v_7v_6v_1v_2v_3v_{10}$ is a $C_8$, $v_{10}\nsim v_4$ otherwise $v_{10}v_9v_8v_4v_{10}$ is a $C_4$, $v_{10}\nsim v_5$ otherwise $v_{10}v_9v_8v_4v_3v_2v_1v_5v_{10}$ is a $C_8$, $v_{10}\nsim v_6$ otherwise $v_{10}v_9v_8v_7v_4v_3v_2v_6v_{10}$ is a $C_8$, $v_{10}\nsim v_7$ otherwise $v_{10}v_9v_8v_7v_{10}$ is a $C_4$. So $v_{10}$ has two neighbors $v_{12},v_{13}\not\in \{v_1,v_2,\ldots,v_{10}\}$. At least one of $v_{12}$ and $v_{13}$, say $v_{12}$, is not adjacent to $v_9$. Moreover, $v_{12}\nsim v_1$ otherwise $v_{12}v_{10}v_{9}v_8v_7v_6v_2v_1v_{12}$ is a $C_8$,  $v_{12}\nsim v_2$ otherwise $v_{12}v_{10}v_{9}v_8v_7v_6v_1v_2v_{12}$ is a $C_8$, $v_{12}\nsim v_5$ otherwise $v_{12}v_{10}v_{9}v_8v_7v_6v_1v_5v_{12}$ is a $C_8$,  $v_{12}\nsim v_6$ otherwise $v_{12}v_{10}v_{9}v_8v_4v_5v_1v_6v_{12}$ is a $C_8$, $v_{12}\nsim v_8$ otherwise $v_{12}v_{10}v_{9}v_8v_{12}$ is a $C_4$. Note that $v_{12}$ can not be adjacent to both $v_{4}$ and $v_7$ since otherwise $v_{12}v_7v_8v_4v_{12}$ is a $C_4$. If $v_{12}\nsim v_7$, then $v_{12}v_{10}v_{9}v_8v_7v_6v_1v_5$ is an induced $P_8$ in $G$, a contradiction. If $v_{12}\nsim v_4$, then $v_{12}v_{10}v_{9}v_8v_4v_5v_1v_2$ is an induced $P_8$ in $G$, a contradiction. (See Figure~1(c) for an illustration.)

{\bf Case 2} $v_7\nsim v_4$.

In this case, $v_7$ has two neighbors $v_8,v_9\not\in\{v_1,v_2,\ldots,v_7\}$. We claim that we may assume $v_8\nsim v_4,v_8\nsim v_6$. Otherwise, if one of $v_8$ and $v_9$, say $v_9$, is adjacent to $v_4$, then $v_8\nsim v_4$ otherwise $v_8v_4v_9v_7v_8$ is a $C_4$ and $v_8\nsim v_6$ otherwise $v_8v_6v_2v_1v_5v_4v_9v_7v_8$ is a $C_8$. So assume that $v_8\nsim v_4,v_9\nsim v_4$.  We can also assume that $v_8\nsim v_6$ since $v_8$ and $v_9$ can not both be adjacent to $v_6$.

Furthermore, $v_8\nsim v_1$ otherwise $v_8v_7v_6v_1v_8$ is a $C_4$, $v_8\nsim v_2$ otherwise $v_8v_7v_6v_2v_8$ is a $C_4$, $v_8\nsim v_3$ otherwise $v_8v_7v_6v_2v_1v_5v_4v_3v_8$ is a $C_8$, $v_8\nsim v_5$ otherwise $v_8v_7v_6v_1v_2v_3v_4v_5v_8$ is a $C_8$. So $v_8$ has two neighbors $v_{10},v_{11}\not\in \{v_1,v_2,\ldots,v_8\}$. At least one of $v_{10}$ and $v_{11}$, say $v_{10}$, is not adjacent to $v_7$. And $v_{10}\nsim v_3$ otherwise $v_{10}v_8v_7v_6v_1v_5v_4v_3v_{10}$ is a $C_8$, $v_{10}\nsim v_4$ otherwise $v_{10}v_8v_7v_6v_1v_2v_3v_4v_{10}$ is a $C_8$, $v_{10}\nsim v_5$ otherwise $v_{10}v_8v_7v_6v_2v_3v_4v_5v_{10}$ is a $C_8$, $v_{10}\nsim v_6$ otherwise $v_{10}v_8v_7v_6v_{10}$ is a $C_4$. Note that $v_{10}$ can not be adjacent to both $v_1$ and $v_2$, otherwise $v_{10}v_1v_6v_2v_{10}$ is a $C_4$. If $v_{10}\nsim v_1$, then $v_{10}v_8v_7v_6v_1v_5v_4v_3$ is an induced $P_8$ in $G$, a contradiction. If $v_{10}\nsim v_2$, then $v_{10}v_8v_7v_6v_2v_3v_4v_5$ is an induced $P_8$ in $G$, a contradiction. (See Figure~1(d) for an illustration.)
\end{proof}

\begin{claim}\label{claim2.2} $k\geq 6$.
\end{claim}

\begin{proof}
Suppose that $k=5$. Since $\delta(G)\geq 3$, $v_1$ has a neighbor $v_6\not\in\{v_1,v_2,\ldots,v_5\}$. By Claim~\ref{claim2.1} and $G$ contains no $C_4$, $v_6\nsim v_i$ for $i\in\{2,3,4,5\}$. So $v_6$ has two neighbors $v_7,v_8\not\in\{v_1,v_2,\ldots,v_5\}$.

We claim that we may assume $v_7\nsim v_3,v_7\nsim v_4$.  If one of $v_7$ and $v_8$, say $v_8$, is adjacent to $v_3$, then $v_7\nsim v_3$ otherwise $v_7v_3v_8v_6v_7$ is a $C_4$, and $v_7\nsim v_4$ otherwise $v_7v_6v_8v_3v_2v_1v_5v_4v_7$ is a $C_8$. By the symmetry between $v_7$ and $v_8$, we then assume that $v_7\nsim v_3$ and $v_8\nsim v_3$. Furthermore, we can assume that $v_7\nsim v_4$ since $v_7$ and $v_8$ can not be both adjacent to $v_4$.

We can also assume  $v_7\nsim v_1$. Otherwise, suppose that $v_7\thicksim v_1$. Then $v_8\nsim v_1$ otherwise $v_7v_6v_8v_1v_7$ is a $C_4$, $v_8\nsim v_3$ otherwise take $C=v_1v_2v_3v_8v_6v_1$ and we obtain a contradiction to Claim~\ref{claim2.1}, $v_8\nsim v_4$ otherwise take $C=v_1v_6v_8v_4v_5v_1$ and we obtain a contradiction to Claim~\ref{claim2.1}. So we take $v_8$  to play the role of $v_7$.

Furthermore, $v_7\nsim v_2,v_7\nsim v_5$, otherwise there is a $C_4$ in $G$. By the  discussion above, $v_7$ has two neighbors $v_9,v_{10}\not\in \{v_1,v_2,\ldots,v_7\}$. We claim that we may assume $v_9\nsim v_i$, $i\in \{1,2,\ldots,6\}$. It is easy to check that $v_9\nsim v_1$ otherwise $v_9v_7v_6v_1v_9$ is a $C_4$, $v_9\nsim v_2$ otherwise $v_9v_7v_6v_1v_5v_4v_3v_2v_9$ is a $C_8$, $v_9\nsim v_5$ otherwise $v_9v_7v_6v_1v_2v_3v_4v_5v_9$ is a $C_8$. By symmetry, $v_{10} \not\sim v_1, v_2, v_5$.
 If one of $v_9$ and $v_{10}$, say $v_{10}$, is adjacent to $v_6$, then $v_9\nsim v_6$ otherwise $v_9v_7v_{10}v_{6}v_9$ is a $C_4$, $v_9\nsim v_3$ otherwise $v_9v_7v_{10}v_6v_1v_5v_4v_3v_9$ is a $C_8$, $v_9\nsim v_4$ otherwise $v_9v_7v_{10}v_6v_1v_2v_3v_4v_9$ is a $C_8$. By the symmetry between $v_9$ and $v_{10}$, we assume that $v_9\nsim v_6$ and $v_{10}\nsim v_6$. If one of $v_9$ and $v_{10}$, say $v_{10}$, is adjacent to $v_3$, then $v_9\nsim v_3$ otherwise $v_9v_7v_{10}v_3v_9$ is a $C_4$, $v_9\nsim v_4$ otherwise $v_9v_7v_{10}v_3v_2v_1v_5v_4v_9$ is a $C_8$. So assume that $v_9\nsim v_3,v_{10}\nsim v_3$. Finally we can assume that $v_9\nsim v_4$ since $v_9$ and $v_{10}$ can not be both adjacent to $v_4$.

\begin{center}
\begin{tikzpicture}

{\tikzstyle{every node}=[draw,circle,fill=white,minimum size=4pt,
                            inner sep=0pt]
    \draw (3,6) node (v1) {}
        -- ++(0:1cm) node (v2){}
        -- ++(0:1cm) node (v3){}
       -- ++(0:1cm) node (v4){}
        -- ++(0:1cm) node (v5){};

 \draw (v1)  {}
 -- ++(270:1cm) node (v6)   {}
  -- ++(240:1cm) node (v7)   {}
 --++(240:1cm) node (v9)   {}
  -- ++(240:1cm) node (v11)   {} ;
 \draw (v6)  {}
 -- ++(300:1cm) node (v8)   {} ;
  \draw (v7)  {}
 -- ++(300:1cm) node (v10)   {} ;
 \draw (v9)  {}
 -- ++(300:1cm) node (v12)   {} ;
 \draw (v1) .. controls (4,7) and (7,6.5)..(v5) ;
        }

 \draw (v1) node [left] {$v_1$};
 \draw (v2) node [below] {$v_2$};
 \draw (v3) node [below] {$v_3$};
 \draw (v4) node [below] {$v_4$};
 \draw (v5) node [below] {$v_5$};
 \draw (v6) node [left] {$v_6$};
 \draw (v7) node [left] {$v_7$};
 \draw (v8) node [right] {$v_8$};
 \draw (v9) node [left] {$v_9$};
 \draw (v10) node [right] {$v_{10}$};
 \draw (v11) node [left] {$v_{11}$};
 \draw (v12) node [right] {$v_{12}$};

\node at (5,1.6) {{\bf Figure 2:} Illustration for Claim~\ref{claim2.2}};
\end{tikzpicture}
\end{center}

By the assumption that $v_9\nsim v_i$, $i\in \{1,2,\ldots,6\}$, $v_9$ has has two neighbors $v_{11},v_{12}\not\in \{v_1,v_2,\ldots,v_7\}$. We claim that by the symmetry between $v_{11}$ and $v_{12}$, we may assume that $v_{11}\nsim v_i$, $i\in \{7,6,1,3,4\}$. It is easy to check that $v_{11}\nsim v_6$ otherwise $v_{11}v_9v_7v_6v_{11}$ is a $C_4$, $v_{11}\nsim v_3$ otherwise $v_{11}v_9v_7v_6v_1v_5v_4v_3v_{11}$ is a $C_8$, $v_{11}\nsim v_4$ otherwise $v_{11}v_9v_7v_6v_1v_2v_3v_4v_{11}$ is a $C_8$. Symmetrically, $v_{12} \not\sim v_6, v_3, v_4$.
If one of $v_{11}$ and $v_{12}$, say $v_{12}$, is adjacent to $v_7$, then $v_{11}\nsim v_1$ otherwise take $C=v_{11}v_9v_7v_6v_1v_{11}$ and we obtain a contradiction to Claim~\ref{claim2.1}, $v_{11}\nsim v_7$ otherwise $v_{11}v_9v_{12}v_7v_{11}$  is a $C_4$. So assume that $v_{11}\nsim v_7,v_{12}\nsim v_7$. Finally we can assume that $v_{11}\nsim v_1$ since $v_{11}$ and $v_{12}$ can not be both adjacent to $v_1$.

Note that $v_{11}\thicksim v_2$ and $v_{11}\thicksim v_5$ can not be both hold since otherwise $v_{11}v_2v_{1}v_5v_{11}$ is a $C_4$. If $v_{11}\nsim v_5$, then $v_{11}v_9v_7v_6v_1v_5v_4v_3$ is an induced $P_8$ in $G$, a contradiction. If $v_{11}\nsim v_2$, then $v_{11}v_9v_7v_6v_1v_2v_3v_4$ is an induced $P_8$ in $G$, a contradiction. (See Figure~2 for an illustration.)
\end{proof}

\begin{claim}\label{claim2.3} $k=7$.
\end{claim}

\begin{proof} Suppose that $k\leq 6$, then $k=6$ by Claim~\ref{claim2.2}. So $G$ contains no $C_5$ since $G$ contains no induced $C_5$ and no $C_4$.

{\bf Case 1} The cycle $C$ has two consecutive vertices that have a common neighbor in $G$.

We assume, without loss of generality, that $v_1$ and $v_2$ have a common neighbor $v_7$. Then $v_7\nsim v_i$ for $i\in\{3,4,5,6\}$ as $G$ contains no $C_4,C_5$ or $C_8$. So $v_7$ has a neighbor $v_8\not\in\{v_1,v_2,\ldots,v_7\}$. Then $v_8\nsim v_i$ for $i\in\{1,2,\ldots,6\}$ as $G$ contains no $C_4,C_5$ or $C_8$. It follows that $v_8$ has two neighbors $v_9,v_{10}\not\in\{v_1,v_2,\ldots,v_8\}$. At least one of $v_9$ and $v_{10}$, say $v_9$, is not adjacent to $v_7$ since there is no $C_4$ in $G$. Moreover, $v_9\nsim v_i$ for $i\in\{1,2,\ldots,6\}$ as $G$ contains no $C_4,C_5$ or $C_8$. So $v_9v_8v_7v_2v_3v_4v_5v_6$ is an induced $P_8$ in $G$, a contradiction. (See Figure~3(a) for an illustration.)

\begin{center}
\begin{tikzpicture}

{\tikzstyle{every node}=[draw,circle,fill=white,minimum size=4pt,
                            inner sep=0pt]
    \draw (1,6) node (v1) {}
        -- ++(0:1cm) node (v2){}
        -- ++(0:1cm) node (v3){}
       -- ++(0:1cm) node (v4){}
        -- ++(0:1cm) node (v5){}
        -- ++(0:1cm) node (v6){};

 \draw (v1)  {}
 -- ++(300:1cm) node (v7)   {}
  -- ++(270:1cm) node (v8)   {}
 --++(240:1cm) node (v9)   {} ;
 \draw (v8)  {}
 -- ++(300:1cm) node (v10)   {} ;
 \draw (v7)--(v2);
 \draw[red] (v9)-- (v8) -- (v7)-- (v2)-- (v3)-- (v4)-- (v5)-- (v6);
 \draw (v1) .. controls (3,7) and (5.5,6.5)..(v6) ;
        }

 \draw (v1) node [left] {$v_1$};
 \draw (v2) node [below right] {$v_2$};
 \draw (v3) node [below] {$v_3$};
 \draw (v4) node [below] {$v_4$};
 \draw (v5) node [below] {$v_5$};
 \draw (v6) node [below] {$v_6$};
 \draw (v7) node [left] {$v_7$};
 \draw (v8) node [left] {$v_8$};
 \draw (v9) node [left] {$v_9$};
 \draw (v10) node [right] {$v_{10}$};

\node at (3.5,2.5) {(a) Case 1};

{\tikzstyle{every node}=[draw,circle,fill=white,minimum size=4pt,
                            inner sep=0pt]
    \draw (8,6) node (v1) {}
        -- ++(0:1cm) node (v2){}
        -- ++(0:1cm) node (v3){}
       -- ++(0:1cm) node (v4){}
        -- ++(0:1cm) node (v5){}
        -- ++(0:1cm) node (v6){};

 \draw (v1)  {}
 -- ++(270:1cm) node (v7)   {}
  -- ++(240:1cm) node (v8)   {}
 --++(240:1cm) node (v10)   {} ;
 \draw (v7)  {}
 -- ++(300:1cm) node (v9)   {} ;
  \draw (v8)  {}
 -- ++(300:1cm) node (v11)   {} ;
 \draw[red] (v10)-- (v8) -- (v7)-- (v1)-- (v2)-- (v3)-- (v4)-- (v5);
 \draw (v1) .. controls (10,7) and (12.5,6.5)..(v6) ;
        }

 \draw (v1) node [left] {$v_1$};
 \draw (v2) node [below right] {$v_2$};
 \draw (v3) node [below] {$v_3$};
 \draw (v4) node [below] {$v_4$};
 \draw (v5) node [below] {$v_5$};
 \draw (v6) node [below] {$v_6$};
 \draw (v7) node [left] {$v_7$};
 \draw (v8) node [left] {$v_8$};
 \draw (v9) node [right] {$v_9$};
 \draw (v10) node [left] {$v_{10}$};
  \draw (v11) node [right] {$v_{11}$};

\node at (10.5,2.5) {(b) Case 2};

\node at (7,1.8) {{\bf Figure 3:} Illustration for Claim~\ref{claim2.3}};
\end{tikzpicture}
\end{center}

{\bf Case 2} No two consecutive vertices on $C$ share a common neighbor in $G$.

Since $\delta(G)\geq 3$, $v_1$ has a neighbor $v_7\not\in\{v_1,v_2,\ldots,v_6\}$. Then $v_7\nsim v_i$ for $i\in\{2,3,\ldots,6\}$ since $G$ contains no $C_4$ or $C_5$ and by the assumption of Case 2. So $v_7$ have two neighbors $v_8,v_9\not\in\{v_1,v_2,\ldots,v_7\}$. We claim that by the symmetry between $v_8$ and $v_9$, we may assume that $v_8\nsim v_1$ and $v_8\nsim v_4$.  If one of $v_8$ and $v_9$, say $v_9$, is adjacent to $v_1$, then $v_8\nsim v_1$ otherwise $v_8v_1v_9v_7v_8$ is a $C_4$, $v_8\nsim v_4$ otherwise we take $C=v_8v_7v_1v_2v_3v_4v_8$ and it is back to Case 1. So we assume $v_8\nsim v_1$ and $v_9\nsim v_1$. We can assume that $v_8\nsim v_4$ since $v_8$ and $v_9$ can not be both adjacent to $v_4$.

Moreover $v_8\nsim v_i$ for $i\in\{2,3,5,6\}$ since $G$ contains no $C_4,C_5$ or $C_8$. So $v_8$ has two neighbors $v_{10},v_{11}\not\in\{v_1,v_2,\ldots,v_8\}$. By the symmetry between $v_{10}$ and $v_{11}$, we  claim that we may assume $v_{10}\nsim v_7,v_{10}\nsim v_4$. If one of $v_{10}$ and $v_{11}$, say $v_{11}$, is adjacent to $v_7$, then $v_{10}\nsim v_7$ otherwise $v_{10}v_8v_{11}v_7v_{10}$ is a $C_4$, $v_{10}\nsim v_4$ otherwise $v_{10}v_8v_{11}v_7v_1v_2v_3v_4v_{10}$ is a $C_8$. So we assume $v_{10}\nsim v_7$ and $v_{11}\nsim v_7$. Finally we may assume that $v_{10}\nsim v_4$ since $v_{10}$ and $v_{11}$ can not be both adjacent to $v_4$.

Furthermore, $v_{10}\nsim v_i$ for $i\in \{1,2,3,5,6\}$ since $G$ contains no $C_4,C_5$ or $C_8$. Then $v_{10}v_8v_7v_1v_2v_3v_4v_5$ is an induced $P_8$, a contradiction. (See Figure~3(b) for an illustration.)
\end{proof}

By Claim~\ref{claim2.3}, $k=7$. Then $G$ contains no $C_i$ for $i\in \{4,5,6,8\}$. Since $\delta(G)\geq 3$, $v_1$ has a neighbor $v_8\not\in\{v_1,v_2,\ldots,v_7\}$. Then $v_8\nsim v_i$ for $i\in\{2,3,\ldots,7\}$. It follows that $v_8$ has two neighbors $v_9,v_{10}\not\in\{v_1,v_2,\ldots,v_8\}$. We assume, without loss of generality, that $v_{9}\nsim v_1$. Since $G$ contains no $C_i$ for $i\in \{4,5,6,8\}$, then $v_9\nsim v_i$ for $i\in\{2,3,\ldots,7\}$. It follows that $v_{9}v_8v_1v_2v_3v_4v_5v_6$ is an induced $P_8$ in $G$, a contradiction. (See Figure~4 for an illustration.)

\begin{center}
\begin{tikzpicture}

{\tikzstyle{every node}=[draw,circle,fill=white,minimum size=4pt,
                            inner sep=0pt]
    \draw (0,6) node (v1) {}
        -- ++(0:1cm) node (v2){}
        -- ++(0:1cm) node (v3){}
       -- ++(0:1cm) node (v4){}
        -- ++(0:1cm) node (v5){}
        -- ++(0:1cm) node (v6){}
        -- ++(0:1cm) node (v7){};

 \draw (v1)  {}
 -- ++(270:1cm) node (v8)   {}
  -- ++(240:1cm) node (v9)   {};
 \draw (v8)  {}
 -- ++(300:1cm) node (v10)   {} ;
 \draw[red] (v9)-- (v8) --  (v1)-- (v2)-- (v3)-- (v4)-- (v5)--(v6);
 \draw (v1) .. controls (1,7) and (5.5,6.5)..(v7) ;
        }

 \draw (v1) node [left] {$v_1$};
 \draw (v2) node [below] {$v_2$};
 \draw (v3) node [below] {$v_3$};
 \draw (v4) node [below] {$v_4$};
 \draw (v5) node [below] {$v_5$};
 \draw (v6) node [below] {$v_6$};
 \draw (v7) node [below] {$v_7$};
 \draw (v8) node [left] {$v_8$};
 \draw (v9) node [left] {$v_9$};
 \draw (v10) node [right] {$v_{10}$};
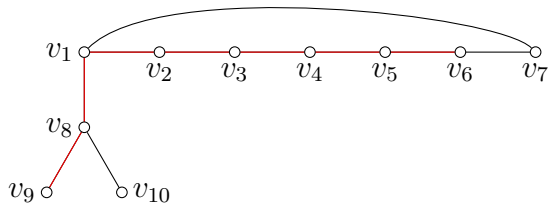
\node at (3,3.5) {{\bf Figure 4:} Illustration for $k=7$};
\end{tikzpicture}
\end{center}

We complete the proof of Theorem~\ref{thm2}.
\end{proof}

%
%

\bibliography{cycle}

\end{document}